\documentclass[12pt]{report}
\usepackage{amssymb,amsmath,mathrsfs}
\usepackage[greek,francais,french]{babel}
\usepackage{geometry}
\usepackage{lscape}
\usepackage{multicol}
\usepackage[T1]{fontenc}
\usepackage[latin1]{inputenc}
\usepackage{xspace}
\usepackage{epsfig,longtable}
\usepackage{url,color,srcltx,mathrsfs}
\usepackage[all]{xy}
   
\usepackage{graphicx}
\usepackage{array}
\usepackage{titlesec}

\makeatletter \@addtoreset{equation}{section}
\@addtoreset{section}{part}

\makeatother

\newtheorem{theorem}{Th\'eor\`eme}[section]

\newtheorem{lemme}[theorem]{Lemme}
\newtheorem{remark}[theorem]{Remarque}

\newtheorem{definition}[theorem]{D\'efinition}

\newenvironment{proof}{\noindent\textbf{Preuve \\ }}{\hspace*{\fill}$\Box$\medskip}

\begin{document}

\title{Sur l'op{\'e}rateur $\bar{\partial}$ et les fonctions
diff{\'e}rentiables au sens de Whitney dans un domaine q-convexe de
$\mathbb{C}^n$}
\author{Eramane Bodian et Salomon Sambou}
\maketitle

\abstract{On r{\'e}sout le probl{\`e}me du
$\bar{\partial}$ dans un domaine $q$-convexe pour
les formes diff{\'e}rentielles dont les co{\'e}fficients sont des fonctions
diff{\'e}rentiables au sens de whitney.
\begin{center}
\textbf{Abstract}
\end{center}                           
We solve the
$\bar{\partial}$-problem on $q$-convex domain for
forms witch coefficients are differentiable functions in the sense of
Whitney.\\
\textbf{Mots cl{\'e}s.}  Domaine $q$-convexe $\bar{\partial}$, fonctions diff{\'e}rentiables au sens de  
Whitney.\\
Classification math{\'e}matique 2010 :
32F32.}

\section{Introduction}

Dans ce papier nous cherchons {\`a} {\'e}tendre les r{\'e}sultats de Dufresnoy
dans \cite{Duf} {\`a} des domaines q-convexes de
$\mathbb{C}^n$. Le passage de la pseudoconvexit{\'e} {\`a} la q-convexit{\'e}
nous oblige {\`a} affaiblir l'hypoth{\`e}se $( \lambda)$ de \cite{Duf} en $( \lambda')$.

Soit $\Gamma$ un ferm{\'e} de $\mathbb{C}^n$ et $W ( \Gamma)$ l'espace des
fonctions infiniment diff{\'e}rentiables au sens de Whitney sur $\Gamma$ , qui
s'identifie au quotient $\frac{\mathcal{C}^{\infty} (
\mathbb{C}^n)}{\mathfrak{F} ( \Gamma)}$ o{\`u} $\mathfrak{F} ( \Gamma)$
d{\'e}signe l'id{\'e}al des fonctions identiquement nulles sur $\Gamma$ ainsi
que toutes leurs d{\'e}riv{\'e}es. $W^{( p, r)} ( \Gamma)$ est l'espace des
formes diff{\'e}rentielles de type $( p, r)$ {\`a} coefficients dans $W (
\Gamma)$. Dans $\mathbf{[ \mathbf{1}]}$, Dufresnoy a montr{\'e} que si
$\Gamma$ poss{\`e}de la propri{\'e}t{\'e} $( \lambda)$ i.e si pour tout $R >
0$ , il existe une suite $( \Omega_{\nu}^R)_{\nu}$ \ d'ouverts pseudo-convexes
de $\mathbb{C}^n$ telle que:
\begin{enumerate}
  \item $\cap \Omega_{\nu}^R = \Gamma_R$=$\{ z \in \Gamma ; | z | \leqslant R
  \}$ ,
  
  \item il existe $p$ ( d{\'e}pendant {\'e}ventuellement de $R$ ) tel que , $0
  \leqslant \varepsilon \leqslant \tfrac{1}{2}$ , il existe $\nu$ avec $\{ z
  \in \mathbb{C}^n ; \rm{dist} ( z, \Gamma) \leqslant \varepsilon^p \}
  \subset \Omega_{\nu}^R \subset \{ z \in \mathbb{C}^n ; \rm{dist} ( z,
  \Gamma) < \varepsilon \}$;
  
  \item si $R' > R$ , pour tout $\mu$ , il existe $\nu_0$ tel que , si $\nu >
  \nu_0$ , $\Omega_{\nu}^R$ \ est holomorphiquement convexe dans
  $\Omega_{\mu}^{R'}$.
\end{enumerate}
Alors pour tout $w \in W^{( p, q)} ( \Gamma)$ , $q \geqslant 1$ avec
$\bar{\partial} w = 0$ , il existe $\alpha \in W^{( p, q - 1)} ( \Gamma)$
telle que $\bar{\partial} \alpha = w$.

Dans notre cas l'hypoth{\`e}se $( \lambda')$ consiste {\`a} remplacer dans $(
\lambda)$ la pseudo-convexit{\'e} par la q-convexit{\'e} et l'hypoth{\`e}se
iii. par la suivante :

iii. Si $R' > R$ , pour tout $\mu$ , il existe $\nu_0$ tel que , si $\nu >
\nu_0$ , les $( 0, q - 1)$-formes diff{\'e}rentielles
$\bar{\partial}$-ferm{\'e}es dans $\overline{\Omega_{\nu}^R}$ sont limites
uniformes de $( 0, q - 1)$-formes diff{\'e}rentielles
$\bar{\partial}$-ferm{\'e}es dans $\Omega_{\mu}^{R'}$.

Dans la suite on consid{\'e}rera les $( 0, r)$-formes diff{\'e}rentielles et
on les notera $r$-formes diff{\'e}rentielles et l'espace $W^{( 0, r)} (
\Gamma)$ par $W^r ( \Gamma)$.

Kohn \cite{Kohn} et L. Ho \cite{Ho} ont obtenu respectivement une r{\'e}solution
globale du $\bar{\partial}$ avec estimation jusqu'au bord pour les domaines
faiblement pseudoconvexes respectivement faiblement $q$-convexes , mais Ho \cite{Ho} a {\'e}tabli ce r{\'e}sultat pour certains classes d'ouverts $q$-convexes.
Malgr{\'e} le fait que ces classes \ ne contiennent pas tous les ouverts
strictement $q$-convexes (remarque de V. Michel \cite{Vin}) , on peut en utilisant
ses ouverts avoir une r{\'e}solution globale dans le cas $q$-convexe du
th{\'e}or{\`e}me principal de Dufresnoy \cite{Duf}. Ainsi on a :

\begin{theorem}\label{th1}
  Soit $\Gamma$ un ferm{\'e} de $\mathbb{C}^n $ poss{\'e}dant la
propri{\'e}t{\'e} $( \lambda^{\prime})$ , alors pour tout $\omega \in W^r ( \Gamma)$
, \ $r \geqslant q$ telle que \ $\bar{\partial} \omega = 0$ , il existe
$\alpha \in W^{( r - 1)} ( \Gamma) $ tel que $\bar{\partial} \alpha = \omega$.
\end{theorem}

Pour avoir un r{\'e}sultat consernant une $q$-convexit{\'e} qui contient tous
les classes d'ouverts strictement $q$-convexes, nous avons utilis{\'e} les
r{\'e}sultats locaux de J.Brinkschulte \cite{Jud1}. On obtient ainsi une version
locale $q$-convexe du th{\'e}or{\`e}me principal dans \cite{Duf} :

\begin{theorem}\label{th2}
  Soit $\Gamma$ un ferm{\'e} de $\mathbb{C}^n$ poss{\'e}dant la
propri{\'e}t{\'e} $( \lambda')$ et $z_0 \in \Gamma$ , alors il existe un
syst{\`e}me fondamental de voisinages $U$ de $z_0$ dans \ $\mathbb{C}^n$ et
pour tout $f \in W^r ( \Gamma \cap \bar{U})$ , $r \geqslant q$ avec
$\bar{\partial} f = 0$ , il existe $\alpha \in W^{r - 1} ( \Gamma \cap
\bar{U})$ telle que $\bar{\partial} \alpha = f$.
\end{theorem}

Notre papier est organis{\'e} comme suit: la section 2 concerne le
r{\'e}sultat global \ avec l'utilisation des r{\'e}sultats de Ho \cite{Ho}, la
section 3 concerne la r{\'e}solution locale en utilisant ceux de J.
Brinkschulte \cite{Jud1} et enfin dans la section 4 on donne quelques exemples de
ferrm{\'e}s de $\mathbb{C}^n$ v{\'e}rifiant la condition $( \lambda')$.

\section{R{\'e}solution globale du $\bar{\partial}$}

\subsection{Pr{\'e}liminaires}

Soit $\Omega$ un domaine de $\mathbb{C}^n$ {\`a} bord lisse et $\rho$ une
fonction de classe $\mathcal{C}^{\infty} $ d{\'e}finie sur $\partial \Omega$
tel que $\rho < 0$ dans $\Omega$ et $| \partial \rho | = 1$ sur le bord. ici
$L_1 , L_2 , \ldots .
, L_n $ sont les coordonn{\'e}es sp{\'e}ciales du bord prises sur un
voisinage $U$ de $x_0 \in \partial \Omega$ i.e $L_i \in T^{1, 0}$ sur $U \cap
\bar{\Omega}$ avec $L_i$ tangentiel pour $1 \leqslant i \leqslant n - 1$ et
$\langle L_i, \overline{L_{}}_j \rangle = \delta_{\rm{ij}}$ . Les duaux
basics des $( 1, 0)$ formes sont $w_1 , w_2, \ldots . , w_n$
avec $w_n = \partial \rho$. Alors $C_{\rm{ij}} = \langle L_i \wedge
\bar{L}_j , \partial \bar{\partial} \rho \rangle , i, j = 1,
\ldots ., n - 1$ est la forme de Levi associ{\'e}e {\`a} $\rho$. On note
$L^2_r ( \Omega)$ \ l'espace des $r$-formes {\`a} coefficients dans $L^2 (
\Omega)$ i.e \ $\left\{ f : \int_{\Omega} | f |^2 d \nu < + \infty \right\}$
avec $f = \sum_{| J | = r} f_J d \bar{z}_J$ une $r$-forme et $d \nu$ la mesure
de Lebesgue. $L^2_r ( \Omega, \rm{loc})$ l'espace des $r$-formes {\`a}
coefficients de carr{\'e} localement int{\'e}grables sur tout compact de
$\Omega$. $L^2_r ( \Omega, \varphi) = \left\{ f : \int_{\Omega} | f |^2 e^{-
\varphi} d \nu < + \infty \right\} $ et on note $\| f \|_{( \varphi)}$ la
norme associ{\'e}e issue du produit scalaire $\langle f, g \rangle =
\int_{\Omega} f \bar{g} e^{- \varphi} d \nu$.$ \overline{\partial^{\ast}}$ est
l'op{\'e}rateur adjoint de $\bar{\partial}$ dans l'espace $L^2_r ( \Omega)$.
$\mathcal{A}_{( r)} ( \Omega)$ est l'espace des $r$-formes diff{\'e}rentielles
dans $\Omega$ de classe $\mathcal{C}^{\infty}$ jusqu'au bord et
$\mathcal{D}_{( r)} ( \Omega)$ l'espace des $r$-formes diff{\'e}rentielles de
classe $\mathcal{C}^{\infty}$ {\`a} support compact dans $\Omega$.

\begin{definition}
  
\end{definition}

\subsubsection{Fonctions q-sousharmoniques ( cf \cite{Ho}.}

Soit $\phi$ une fonction de classe $\mathcal{C}^2$ d{\'e}finie sur un domaine
$U \subseteqq \mathbb{C}^n$. Pour chaque $q \geqslant 1$ on d{\'e}finit une
matrice carr{\'e}e $\Phi^{( q)} ( x)$ d'ordre $\frac{n!}{q! ( n - q) !}$
associ{\'e}e {\`a} $\phi$. En fait les valeurs de la matrice sont les $\phi_{I
J}$ o{\`u} $I$ et $J$ sont des $q$-uplets croissants d'entier entre $1$ et
$n$. On d{\'e}finit
\[ \phi_{I J} ( x) = \left\{ \begin{array}{l}
     \sum_{i \in I} \frac{\partial^2 \phi}{\partial z_i \partial \bar{z}_j} 
     \text{\rm{si}} I = J,\\
     \varepsilon_{i K}^I \varepsilon_{j K}^J \frac{\partial^2 \phi}{\partial
     z_i \partial \bar{z}_j}  \text{\rm{si}} I = \langle i K \rangle
     , \text{} J = \langle j K \rangle, \text{\rm{et}} i \neq j,\\
     0 \text{\rm{sinon}},
   \end{array} \right. \]
o{\`u} $\langle i K \rangle$ est l'ordre croissant des indices de l'ensemble
$\{ i \} \cup K$ , et $\varepsilon_{i K}^I$ est le signe de la permutation de
$i K$ {\`a} $I$ , qui est {\'e}gal {\`a} $0$ si $\langle i K \rangle \neq I$.

On dit que $\phi$ est $q$-sousharmonique respectivement strictement
$q$-sousharmonique sur $U \subseteqq \mathbb{C}^n$ si la matrice associ{\'e}e
$\Phi^{( q)} ( x)$ est semi-d{\'e}finie positive respectivement d{\'e}finie
positive pour tout $x \in U$.

\subsubsection{Domaine $q$-convexe (au sens de Ho cf \cite{Ho} ).}

Soit $\Omega$ un domaine {\`a} bord lisse de $\mathbb{C}^n$ et $\rho$ une
fonction d{\'e}finissante de $\Omega$ , alors on dit que $\Omega$ est
$q$-convexe si en chaque point $x_0 \in \partial \Omega$ on a
\[ \sum_K' \sum_{i, j = 1}^n \frac{\partial^2 \rho}{\partial z_i \partial
   \bar{z}_j} u_{\rm{iK}} \bar{u}_{\rm{jK}} \geqq 0 \text{\rm{pour}
   \rm{toute}} q \text{- \rm{forme}} u = \sum_{| K | = q}' u_K d \bar{z}_k
   \text{\rm{telle} \rm{que}} \sum^n_{i = 1} \frac{\partial
   \rho}{\partial z_i} u_{\rm{iK}} = 0 \]
pour tout $| K | = q - 1$.

Le r{\'e}sultat ci-dessus est une adaptation de celui de Dufresnoy \cite{Duf}
(lemme 1):

\begin{lemme}\label{l1}
  
\end{lemme}

Soit $V$ un ouvert q-convexe born{\'e} de $\mathbb{C}^n$ et d{\'e}signons
pour tout $\varepsilon > 0$ par $V^{\varepsilon} = \{ z \in V ; d ( z, V^c) >
\varepsilon \}$ o{\`u} $V^c$ d{\'e}signe le compl{\'e}mentaire de $V$ dans
$\mathbb{C}^n$. Pour toute $r$-forme $f$ {\`a} coefficients dans
$\mathcal{C}^{\infty} ( \bar{V})$ telle que $\bar{\partial} f = 0$ , il existe
$u \in \mathcal{C}_{r - 1}^{\infty} ( V)$ telle que , pour tout $s \in
\mathbb{N}$ et tout $\varepsilon > 0$, $\| u \|_{( s + 1, V^{\varepsilon})}
\leqslant \frac{M_s}{\varepsilon^{s + 1}} \| f \|_{( s, V)}$ o{\`u} $M_s$ ne
d{\'e}pend que du diam{\`e}tre de $V$.\\

\noindent Pr{\'e}cisons que
\[ \| u \|_{( s, V)}^2 = \sum_{| \alpha | \leqslant s} \int_V | D^{\alpha} u
   |^2 d \nu . \]
   
Pour faire la preuve de ce lemme \ref{l1}, on a besoin des deux r{\'e}sultats
suivant:

\begin{lemme}\label{l2}
  (cf lemme 2 de \cite{Duf})
\end{lemme}

Soit $F_1$ et $F_2$ deux ferm{\'e}s de $\mathbb{R}^n$ tels que $d ( F_1, F_2)
\geqslant \delta$. Alors il existe $\varphi \in \mathcal{C}^{\infty} (
\mathbb{R}^n)$ telle que $\varphi$ soit {\'e}gale {\`a} $1$ au voisinage de
$F_1$, $\varphi$ soit nulle au voisinage de $F_2$ et v{\'e}rifie de plus, pour
tout multiindice $\alpha$
\[ \sup_{x \in \mathbb{R}^n} | D^{\alpha} \varphi ( x) | \leqslant
   \frac{N_{| \alpha |}}{\delta^{| \alpha |}} \]
(o{\`u} $N_{| \alpha |}$ ne d{\'e}pend pas de $F_1$ et $F_2$).

\noindent Ce r{\'e}sultat est une adaptation d'une version tr{\`e}s particuli{\`e}re du

lemme 4.4.1 de \cite{Horm}, compl{\'e}t{\'e} par la remarque de \cite{Horm} page 87.

\begin{theorem}\label{th3}
  
\end{theorem}

Soit $\varphi = | z |^2$ et $U$ un ouvert $q$-convexe de $\mathbb{C}^n$ ; si
$f$ est une $r$-forme diff{\'e}rentielle {\`a} co{\'e}fficients dans $L^2 (
U)$ telle que $\bar{\partial} f = 0$ , il existe $\omega$ une $( r - 1)$-forme
diff{\'e}rentielle telle que \ $\bar{\partial} \omega = f$ ; \
$\overline{\partial^{\ast}} ( e^{- \varphi} \omega) = 0$ ; \ $\| \omega
\|_{\varphi} \leqslant \| f \|_{\varphi}$.

{\textbf{Preuve Lemme \ref{l1}}}

Nous pouvons supposer, puisque $\bar{\partial}$ est {\`a} coefficients
constants, que $0 \in U$ et on notera dans la suite $\| \mathbf{g} \|_{( s,
\varepsilon)}$ (resp. $\| \mathbf{g} \|_{( s)}$) au lieu de $\| \mathbf{g}
\|_{( s, U^{\varepsilon})}$ (resp. $\| \mathbf{g} \|_{( s, U)}$).
\begin{enumerate}
  \item {\textrm{{\textbf{Le cas $s = 0$}}.}}
  
  On d{\'e}signe par $\chi_{\varepsilon}$ une fonction {\'e}gale {\`a} $1$ au
  voisinage de $U^{\varepsilon}$ et {\`a} support dans $U^{\frac{\varepsilon}{6}}$ fournie par le
  lemme \ref{l2} ; et on d{\'e}signe par $u$ la solution de l'{\'e}quation
  $\bar{\partial} u = f$ fournie par le th{\'e}or{\`e}me \ref{th3}.
  
  En remarquant que $\| \chi_{\varepsilon} u \|_{( 1)} \geqslant \| u \|_{( 1,
  \varepsilon)}$ et que
  \[ \| \chi_{\varepsilon} u \|^2_{( 1)} = \| \chi_{\varepsilon} u \|^2 + \|
     \bar{\partial} \chi_{\varepsilon} u \|^2 + \| \bar{\partial}^{\ast}
     \chi_{\varepsilon} u \|^2 \]
  il suffit de montrer que $\| \chi_{\varepsilon} u \|^2 + \| \bar{\partial}
  \chi_{\varepsilon} u \|^2 + \| \bar{\partial}^{\ast} \chi_{\varepsilon} u
  \|^2 \leqslant \frac{M_0^2}{\varepsilon^2} \| f \|^2$. Il existe des
  constantes $\mu$ et $\nu$, ne d{\'e}pendant que du diam{\`e}tre de $U$
  telles que
  \[ \| \chi_{\varepsilon} u \|^2 \leqslant \| u \|^2 \leqslant \mu \| u
     \|_{\varphi}^2 \leqslant \mu \| f \|_{\varphi}^2 \leqslant \nu \| f \|^2
  \]
  donc
  \[ \| \chi_{\varepsilon} u \|^2 \leqslant \| u \|^2 \leqslant \nu \| f \|^2
  \]
  d'autre part,
  \[ \| \bar{\partial} \chi_{\varepsilon} u \|^2 \leqslant 2 \{ \|
     \bar{\partial} \chi_{\varepsilon} \wedge u \|^2 + \| \chi_{\varepsilon}
     \bar{\partial} u \|^2 \} \leqslant 2 \left\{ \frac{N_1^2}{\varepsilon^2}
     \| u \|^2 + \| f \|^2 \right\} . \]
  soit encore
  \[ \| \bar{\partial} \chi_{\varepsilon} u \|^2 \leqslant
     \frac{K_1^2}{\varepsilon^2} \| f \|^2 . \]
  Enfin, on a
  \[ \bar{\partial}^{\ast} \chi_{\varepsilon} u = \chi_{\varepsilon}
     \bar{\partial}^{\ast} \chi_{\varepsilon} u + [ \bar{\partial}^{\ast},
     \chi_{\varepsilon}] u. \]
  La condition $\bar{\partial}^{\ast} ( e^{- \varphi} u) = 0$ se traduit par
  le fait que $\bar{\partial}^{\ast}$ agit sur $u$ comme un op{\'e}rateur
  d'ordre $0$ dont les coefficients sont major{\'e}s par une constante ne
  d{\'e}pendant que du diam{\`e}tre de $U$ et $\frac{L_1}{\varepsilon} .$ On a
  donc :
  \[ \| \bar{\partial}^{\ast} ( \chi_{\varepsilon} u) \|^2 \leqslant 2
     \left\{ L_1'^2 \| u \|^2 + \frac{L_1^2}{\varepsilon^2} \| u \|^2 \right\}
     \leqslant 2 \left\{ \nu L_1'^2 \| f \|^2 + \nu
     \frac{L_1^2}{\varepsilon^2} \| f \|^2 \right\} . \]
  En ajoutant les trois termes, on obtient le r{\'e}sultat d{\'e}sir{\'e}.
  
  \item {\textrm{{\textbf{Le cas g{\'e}n{\'e}ral.}}}}
  
  Nous avons montr{\'e} que $\| \chi_{\varepsilon} u \|_{( 1)} \leqslant
  \frac{M_0}{\varepsilon} \| f \|$.
  
  Supposons qu'on ait d{\'e}montr{\'e} $\| \chi_{\varepsilon} u \|_{( s + 1)}
  \leqslant \frac{M_s}{\varepsilon^{s + 1}} \| f \|_{( s)}$ et montrons qu'on
  peut en d{\'e}duire le m$\widehat{\text{e}}$me r{\'e}sultat pour $s + 1$,
  {\`a} savoir $\| \chi_{\varepsilon} u \|_{( s + 2)} \leqslant \frac{M_{s +
  1}}{\varepsilon^{s + 2}} \| f \|_{( s + 1)}$. Soit donc $D^{\alpha}$ une
  d{\'e}rivation d'ordre $s + 1$; on a :
  \[ \| D^{\alpha} \chi_{\varepsilon} u \| \leqslant \frac{M_s}{\varepsilon^{s
     + 1}} \| f \|_{( s)} \leqslant \frac{M_s}{\varepsilon^{s + 2}} \| f \|_{(
     s + 2)}  \text{} ( \varepsilon < 1) ; \]
  d'autre part, comme dans le cas $s = 0$ ,
  \[ \| D^{\alpha} \chi_{\varepsilon} u \|_{( 1)}^2 = \| D^{\alpha}
     \chi_{\varepsilon} u \|^2 + \| \bar{\partial} D^{\alpha}
     \chi_{\varepsilon} u \|^2 + \| \bar{\partial}^{\ast} D^{\alpha}
     \chi_{\varepsilon} u \|^2 . \]
  On a
  \[ \bar{\partial} D^{\alpha} \chi_{\varepsilon} u = D^{\alpha}
     \chi_{\varepsilon} \bar{\partial} u + D^{\alpha} ( \bar{\partial}
     \chi_{\varepsilon} \wedge u) . \]
  Pour le premier terme , on a
  \[ \| D^{\alpha} \chi_{\varepsilon} \bar{\partial} u \| \leqslant \|
     \chi_{\varepsilon} f \|_{( s + 1)} \leqslant \frac{T_{s +
     1}}{\varepsilon^{s + 1}} \| f \|_{( s + 1)} . \]
  Pour le deuxi{\`e}me terme , gr$\widehat{\textrm{a}}$ce {\`a} la formule de
  Leibnitz, on a
  \[ D^{\alpha} ( \bar{\partial} \chi_{\varepsilon} \wedge u) = \sum
     C_{\alpha}^{\beta} ( D^{\beta} \bar{\partial} \chi_{\varepsilon}) \wedge
     D^{\alpha - \beta} u \]
  et
  \[ \| ( D^{\beta} \bar{\partial} \chi_{\varepsilon}) \wedge D^{\alpha -
     \beta} u \| \leqslant \frac{N_{| \beta |}}{\varepsilon^{| \beta | + 1}}
     \| D^{\alpha - \beta} u \|_{\left( 0,  \right)} . \]
  Il suffit d'utiliser la r{\'e}currence pour avoir
  \[ \| D^{\alpha - \beta} u \|_{\left( 0,  \right)} \leqslant \frac{M_{|
     \alpha - \beta |} 6^{| \alpha - \beta |}}{\varepsilon^{| \alpha - \beta
     |}} \| f \|_{( | \alpha - \beta | - 1)} . \]
  Apr{\`e}s sommation sur les multiindices $\beta \leqslant \alpha$, on
  obtient
  \[ \| D^{\alpha} ( \bar{\partial} \chi_{\varepsilon} \wedge u) \| \leqslant
     \frac{S_{s + 1}}{\varepsilon^{| \beta | + 1 + | \alpha - \beta |}} \| f
     \|_{( | \alpha |)} = \frac{S_{s + 1}}{\varepsilon^{s + 2}} \| f \|_{( s +
     1)} . \]
  Par analogie pour $\bar{\partial}^{\ast} D^{\alpha} \chi_{\varepsilon} u$,
  on obtient le r{\'e}sultat cherch{\'e}. \ \ \
  
  \ \ \ \ \ \ \ \ \ \ \ \ \ \ \ \ \ \ \ \ \ \ \ \ \ \ \ \ \ \ \ \ \ \ \ \ \
  \ \ \ \ \ \ \ \ \ \ \ \ \ \ \ \ \ \ \ \ \ \ \ \ \ \ \ \ \ \ \ \ \ \ \ \ \ \
  \ \ \ \ \ $\Box$
\end{enumerate}

Pour faire la preuve de notre th{\'e}or{\`e}me \ref{th1} , on a besoin d'un r{\'e}sultat
de r{\'e}solution du $\bar{\partial}$ avec estimation jusqu'au bord. Ainsi Ho
\cite{Ho} a fait la r{\'e}solution dans un anneau et le r{\'e}sultat reste vrai
dans un domaine $q$-convexe de $\mathbb{C}^n$.

\begin{theorem}\label{th4}
  
\end{theorem}

Soit $\Omega_1$ et $\Omega_2$ deux domaines born{\'e}s de $\mathbb{C}^n$ \ de
classe $\mathcal{C}^{\infty}$ avec $\bar{\Omega}_2 \subseteq \Omega_1$ tels
que $\Omega_1$ est $( q - 1)$-convexe et $\Omega_2$ est $( n - q -
1)$-convexe. Soit $\varrho$ une fonction lisse sur $\bar{\Omega}$ telle que
$\varrho = | z |^2$ dans un voisinage de $\partial \Omega_1$ et $- | z |^2$
dans un voisinage de $\partial \Omega_2$. On pose $\Omega = \Omega_1 -
\Omega_2$ , alors pour $n \geqslant 3$ et $1 \leqslant q \leqslant n - 1$ , si
$\alpha \in \mathcal{A}_{( r)} ( \Omega)$ tel que $\bar{\partial} \alpha = 0$
et $\langle \alpha, \psi \rangle = 0$ pour tout $\psi \in \mathcal{D}_{( r)} (
\Omega)$ et $\bar{\partial}^{\ast} \psi = 0$ , alors il existe $u \in
\mathcal{A}_{( r - 1)} ( \Omega)$ tel que $\bar{\partial} u = \alpha$.

Pour se mettre dans le cas q-convexe , il suffit de prendre $\Omega_2 =
\varnothing$.

Ainsi nous pouvons {\'e}tablir la preuve de notre r{\'e}sultat.

\subsection{Preuve du th{\'e}or{\`e}me \ref{th1}}

La preuve est identique {\`a} celle de Dufresnoy \cite{Duf}. En effet faisons les
d{\'e}tails :

\subsubsection{Le cas o{\`u} $\Gamma$ est born{\'e}.}

On choisit un nombre $R$ tel que $\Gamma = \Gamma_R$ et on notera
$\Omega_{\nu}^R = \Omega_{\nu}$. Quitte {\`a} extraire une sous-suite de la
suite $( \Omega_{\nu})_{\nu}$ initiale , on peut supposer qu'il existe $0 <
\eta < \tfrac{1}{2}$ tel que
\[ \{ z \in \mathbb{C}^n ; d ( z, \Gamma) < \eta^{p^{\nu + 1}} \} \subset
   \Omega_{\nu} \subset \{ z \in \mathbb{C}^n ; d ( z, \Gamma) <
   \eta^{p^{\nu}} \} . \]
Soit $\omega \in W^r ( \Gamma)$ , $r \geqslant q$ avec $\bar{\partial} \omega
= 0$. D{\'e}signons par $\tilde{\omega}$ un prolongement de $\omega$ qui soit
{\`a} coefficients dans $\mathcal{C}^{\infty} ( \mathbb{C}^n)$. Il n'y a
aucune raison pour que $\bar{\partial} \tilde{\omega} = 0$ , mais
n{\'e}anmoins, il existe une constante $C_{N, s}$ telle que $\| \bar{\partial}
\tilde{\omega} \|_{( s, \Omega_{\nu})} \leqslant C_{N, s} \eta^{N_{p^{\nu}}}$
(o{\`u} $N_{p^{\nu}}$ est un r{\'e}el d{\'e}pendant de $R$) pour tout entier
$s$ car $\bar{\partial} \tilde{\omega}$ est identiquement nul sur $\Gamma$
ainsi que toutes ses d{\'e}riv{\'e}es. On d{\'e}signe par $h_{\nu}$ une
solution de $\bar{\partial} h_{\nu} = \bar{\partial} \tilde{\omega}$ dans
$\Omega_{\nu}$ fournie par le lemme \ref{l1} ; on a donc
\[ \| h_{\nu} \|_{( s + 1, \Omega_{\nu + 1})} \leqslant M_s \frac{1}{2}
   \eta^{\nu + 1 - ( s + 1)} \| \bar{\partial} \tilde{\omega} \|_{( s,
   \Omega_{\nu})} \]
en remarquant que
\[ \Omega_{\nu + 1} \subset \left\{ z \in \Omega_{\nu} : d ( z,
   \Omega_{\nu}^c) > \frac{1}{2} \eta^{p^{\nu + 1}} \right\} . \]
Consid{\`e}rons alors sur $\Omega_2$ la forme diff{\'e}rentielle
$\tilde{\omega} - h_1$ , on a {\'e}videmment $\bar{\partial} ( \tilde{\omega}
- h_1) = 0$ , il existe donc une solution $\alpha_1$ une $( r - 1)$-forme
diff{\'e}rentielle fournie par le lemme \ref{l1} telle que $\bar{\partial} \alpha_1 =
\tilde{\omega} - h_1$. En outre sur $\Omega_{\nu + 2}$ pour tout $\nu
\geqslant 1$ , on consid{\`e}re la forme diff{\'e}rentielle $h_{\nu} - h_{\nu
+ 1}$ qui est $\bar{\partial}$-ferm{\'e}e et en vertu du lemme \ref{l1} on
d{\'e}signe par $\alpha_{\nu + 1}$ une solution de l'{\'e}quation
$\bar{\partial} \alpha_{\nu + 1} = h_{\nu} - h_{\nu + 1}$. On a donc
\[ \| \alpha_{\nu + 1} \|_{( s + 2, \Omega_{\nu + 3})} \leqslant M_{s + 1}
   \left( \frac{1}{2} \eta^{p^{\nu + 2}} \right)^{- ( s + 2)} M_s \left(
   \frac{1}{2} \eta^{p^{\nu + 1}} \right)^{- ( s + 1)} \| \bar{\partial}
   \tilde{\omega} \|_{( s, \Omega_{\nu})} \]
par suite
\[ \| \alpha_{\nu + 1} \|_{( s + 2, \Omega_{\nu + 3})} \leqslant M_s M_{s +
   1} \left( \frac{1}{2} \eta^{p^{\nu + 2}} \right)^{- ( 2 s + 3)} \|
   \bar{\partial} \tilde{\omega} \|_{( s, \Omega_{\nu})} . \]
Ainsi en majorant $\| \bar{\partial} \tilde{\omega} \|_{( s, \Omega_{\nu})}$ ,
on a
\[ \| \alpha_{\nu + 1} \|_{( s + 2, \Omega_{\nu + 3})} \leqslant M_s M_{s +
   1} C_{N, s} \left( \frac{1}{2} \eta^{p^{\nu + 2}} \right)^{- ( 2 s + 3)}
   \times \eta^{N_{p^{\nu}}} . \]
On choisit alors $N$ suffisamment grand pour que cette in{\'e}galit{\'e} ,
jointe au lemme de Sobolev fournit la convergence de la s{\'e}rie $\sum
\alpha_{\nu}$ dans $W^r ( \Gamma)$ muni de la topologie quotient de celle de
$\mathcal{C}^{\infty} ( \mathbb{R}^n)$.

\subsubsection{Le cas $r > q$.}

Soit $\omega \in W^r ( \Gamma)$ telle que $\bar{\partial} \omega = 0$. Il
existe une suite $( u_n)_n$ avec $u_n \in W^{( r - 1)} ( \Gamma_n)$ , telle
que $\bar{\partial} u_n = \omega_{| \Gamma_n }$ d'apr{\`e}s le
paragraphe 2.2.1. On peut modifier la suite $( u_n)_n$ de telle sorte que
$u_{n + 1 | \Gamma_n } = u_n$ pour tout $n$. Pour s'en convaincre il
suffit de montrer que $u_{n + 1 | \Gamma_n } - u_n$ peut se
prolonger en une $( r - 1)$-forme diff{\'e}rentielle
$\bar{\partial}$-ferm{\'e}e sur $\Gamma_{n + 1}$. En effet , puisque
$\bar{\partial} ( u_{n + 1 | \Gamma_n } - u_n) = 0$ , il existe donc
$v \in W^{( r - 2)} ( \Gamma_n)$ telle que $\bar{\partial} v = u_{n + 1 |
\Gamma_n } - u_n$. Notons $\tilde{v}$ un prolongement de $v$ dans
$\Gamma_{n + 1}$ et $\bar{\partial} \tilde{v}$ est le prolongement
cherch{\'e}.

La s{\'e}rie $\sum u_n$ ainsi modifi{\'e}e converge {\'e}videmment dans $W^{(
r - 1)} ( \Gamma)$ et la somme $u$ de la s{\'e}rie v{\'e}rifie $\bar{\partial}
u = \omega$.

\subsubsection{Le cas $r = q$.}

Pour ce cas , il suffit comme dans \cite{Duf} d'{\'e}tablir le lemme suivant :

\begin{lemme}\label{l3}
  
\end{lemme}

Soit $\Gamma$ un ferm{\'e} de $\mathbb{C}^n $ poss{\'e}dant la
propri{\'e}t{\'e} $( \lambda')$ et deux nombres r{\'e}els $R$ et $R'$ avec $R
< R'$ ; si $f \in W^{q - 1} ( \Gamma_R)$ est telle que $\bar{\partial} f = 0$
alors $f$ est limite dans $W^{q - 1} ( \Gamma_R)$ de $( q - 1)$-formes
diff{\'e}rentielles $\bar{\partial}$-ferm{\'e}es dans $\Omega_1^{R'}$.

\begin{proof}
  Soit \ $\tilde{f}$ une extension de $f$ {\`a} $\mathbb{C}^n$. On a
  $\bar{\partial} \tilde{f}$ est nulle sur $\Gamma_R$ ainsi que toutes ses
  d{\'e}riv{\'e}es. En notant $h_{\nu}$ la solution de l'{\'e}quation
  $\bar{\partial} h_{\nu} = \bar{\partial} \tilde{f}$ dans $\Omega_{\nu}^R$
  fournie par le lemme \ref{l1} , et le fait que $\bar{\partial} \tilde{f}$ est nulle
  sur $\Gamma_R$ ainsi que toutes ses d{\'e}riv{\'e}es, il existe donc une
  constante $C_{N, s}$ telle que $\| \bar{\partial} \tilde{f} \|_{( s,
  \Omega_{\nu})} \leqslant C_{N, s} \eta^{N_{p^{\nu}}}$. Posons $h_0 =
  \tilde{f}$ et consid{\`e}rons la forme diff{\'e}rentielle $h_{\nu} - h_{\nu
  + 1}$ qui est $\overline{\partial}$-ferm{\'e}e sur $\Omega_{\nu + 1}^R$. Il
  existe alors $\alpha_{\nu + 1}$ solution de l'{\'e}quation $\bar{\partial}
  \alpha_{\nu + 1} = h_{\nu} - h_{\nu + 1}$ telle que
  \[ \| \alpha_{\nu + 1} \|_{( s + 2, \Omega_{\nu + 2}^R)} \leqslant M_s M_{s
     + 1} C_{N, s} \left( \frac{1}{2} \eta^{p^{\nu + 2}} \right)^{- ( 2 s +
     3)} \times \eta^{N_{p^{\nu}}} . \]
  On obtient la convergence de la s{\'e}rie $\sum ( h_{\nu} - h_{\nu + 1})_{|
  \Gamma_R }$ dans $W^{( q - 1)} ( \Gamma_R)$ munie de la topologie
  quotient de celle de $\mathcal{C}^{\infty}_r ( \mathbb{C}^n)$. On a ainsi
  \[ f = h_0 - h_{\nu | \Gamma_R } + \sum_{\nu = 1}^{\infty} h_{\nu}
     - h_{\nu + 1 | \Gamma_R } . \]
  Pour $p$ assez grand ,
  \[ f_N = h_0 - h_1 + \sum_{\nu = 1}^N h_{\nu} - h_{\nu + 1 | \Gamma_R
     } \]
  est une $( q - 1)$ forme $\bar{\partial}$-ferm{\'e}e dans $\Omega_p^R$ et
  $f_{N | \Gamma_R } \longrightarrow f$ uniform{\'e}ment dans $W^{q
  - 1} ( \Gamma_R)$.

\end{proof}

\section{R{\'e}solution locale du $\bar{\partial}$}

Notons que dans cette partie on d{\'e}finit une $q$-convexit{\'e} qui contient
toutes les classes d'ouverts strictement $q$-convexes et contient
celle de Ho \cite{Ho}.

\subsection{Pr{\'e}liminaires}

Soit $\Omega \subset \mathbb{C}^n$ un \ ouvert, on \ consid{\`e}re la
fonction continue $\delta = \delta_{\Omega} : \mathbb{C}^n \rightarrow
\mathbb{R}$ d{\'e}finie par :
\begin{eqnarray*}
  \delta ( z) = \delta_{\Omega} ( z) = : \left\{ \begin{array}{l}
    \text{- -} \rm{dist} ( z, \partial \Omega)  \textrm{pour} z \in
    \bar{\Omega},\\
    + \rm{dist} ( z, \partial \Omega)  \textrm{pour} z \notin
    \bar{\Omega}
  \end{array} \right. &  & 
\end{eqnarray*}
O{\`u} dist(.,.) est la distance euclidienne.

Si $\partial \Omega$ est de classe $\mathcal{C}^2$ alors il existe un ouvert
$U$ de $\partial \Omega$ tel que $\delta$ est de classe $\mathcal{C}^2$ dans
$U$.

Soit $f$ une fonction de classe $\mathcal{C}^2$. On note par \ $\lambda_1^f (
z) \leqslant \ldots \leqslant \lambda_n^f ( z)$ les valeurs propres de la
forme de L{\'e}vi
\[ \mathcal{L} ( f, z) ( \xi) = \sum_{i, j = 1}^n \frac{\partial^2
   f}{\partial z_i \partial \bar{z}_j} ( z) \xi_i \overline{\xi_j} . \]

\begin{definition}
  
\end{definition}

$\bullet$ Soit $\Omega \subset \mathbb{C}^n$ un domaine {\`a} bord lisse de
classe $\mathcal{C}^2$ et soit $z_0 \in \partial \Omega$. Supposons que $U$
est un voisinage ouvert de $z_0$ , et soit $\varrho \in \mathcal{C}^2 ( U,
\mathbb{R})$ une fonction telle que $\Omega \cap U = \{ \varrho < 0 \}$ , et
$d \varrho \neq 0$ sur $\partial \Omega$. Soit $q \geqslant 1$ un entier , on
dit que $\Omega$ est $q$-convexe (resp. strictement $q$-convexe ) en $z_0$ si la
forme de L{\'e}vi
\[ \mathcal{L} ( \varrho, z) ( \xi) = \sum_{i, j = 1}^n \frac{\partial^2
   \varrho}{\partial z_i \partial \bar{z}_j} ( z) \xi_i \overline{\xi_j}
   , \text{} \xi \in \mathbb{C}^n, \]
admet au moins $( n - q)$ valeurs propres positives (resp. strictement
positives) sur l'espace tangent holomorphe
\[ T_z^{1, 0} ( \partial \Omega) = \left\{ \xi \in \mathbb{C}^n  |
     \sum_{j = 1}^n \frac{\partial \varrho}{\partial z_j} ( z) \xi_j
   = 0 \right\} , \text{} z \in \partial \Omega \]
pour tout $z \in V \cap \partial \Omega$ , o{\`u} $V$ est un voisinage ouvert
de $z_0$ dans $\mathbb{C}^n$. $\Omega$ est $q$-convexe s'il est $q$-convexe
en tout point de $\partial \Omega$. Si $\Omega$ est non born{\'e}, alors il
est dit $q$-complet (au sens de Andreotti et Grauert) si et seulement si
$\Omega$ admet une fonction d'exhaustion de classe $\mathcal{C}^2$ dont sa
forme de Levi admet au moins $( n - q + 1)$ valeurs propres strictement
positives en chaque point de $\Omega$.

Notons que $q$-complet au sens de Andreotti et Grauert pour un domaine
born{\'e} {\`a} bord de classe $\mathcal{C}^2$ est $q$-convexe (voir \cite{Jud1}
page 144).

\begin{remark}
  
\end{remark}

Soit $\Omega \subset \mathbb{C}^n$ un domaine q-convexe au voisinage de $z_0$
et {\`a} bord lisse de classe  $\mathcal{C}^2$. Alors il existe un
voisinage de $z_0$ dans $\bar{\Omega}$ qui admet une base de voisinages de
domaines strictement $q$-convexes (voir lemme 2.1 dans \cite{Jud1}).

\begin{lemme}\label{l4}
  (\cite{Jud1} lemme 2.1)
\end{lemme}

Soit $\Omega \subset \mathbb{C}^n$ un domaine $q$-convexe au voisinage de
$z_0$ et {\`a} bord lisse de classe $\mathcal{C}^2$ alors il existe un ouvert
$U \subset \Omega$ arbitrairement petit pour lequel il existe un voisinage $V$
de $z_0$ dans $\mathbb{C}^n$ et une base de voisinages d'ouverts $(
U_{\varepsilon})_{0 \leqslant \varepsilon \leqslant \varepsilon_0}$ ,
$\varepsilon_0 > 0$ , strictement $q$-convexes avec une fonction
d{\'e}finissante $\rho_{\varepsilon}$ de classe $\mathcal{C}^2$ dans un
voisinage de $\bar{U}$ satisfaisant :
\begin{enumerate}
  \item $U \cap V = \Omega \cap V$ et $U = U_0$ .
  
  \item Il existe $a > 0$ et $b > 0$ ind{\'e}pendants de $\varepsilon$ tels
  que pour $0 < \varepsilon \leqslant \varepsilon_0$ , on a :
  \[ \{ \rm{dist} (  ., U) \leqslant a \varepsilon \} \subset
     \subset U_{\varepsilon} \subset \subset \{ \rm{dist} ( ., U) \leqslant
     b \varepsilon \} . \]
  \item $( \varepsilon, z) \mapsto \rho_{\varepsilon}$ est de classe
  $\mathcal{C}^2$ dans un voisinage de $[ 0 , \varepsilon_0] \times
  \overline{U}$.
  
  \item Il existe des r{\'e}els positifs $a_1$ et $b_1$ ind{\'e}pendants de
  $\varepsilon$ tels que pour $\varepsilon \in [ 0, \varepsilon_0]$ et $z \in
  U_{\varepsilon}$ on a :
  \[ a_1 | \rho_{\varepsilon} ( z) | \leqslant | \delta_{U_{\varepsilon}} ( z)
     | \leqslant b_1 | \rho_{\varepsilon} ( z) | . \]
  \item Il existe $\gamma > 0$ ind{\'e}pendant de $\varepsilon$ tel que pour
  tout $z \in U_{\varepsilon}$ , on a : $\lambda_q^{\rho_{\varepsilon}} ( z)
  \geqslant \gamma^{\varepsilon}$ o{\`u} $\lambda_1^{\rho_{\varepsilon}} ( z)
  \leqslant \lambda_2^{\rho_{\varepsilon}} ( z) \leqslant \ldots \leqslant
  \lambda_q^{\rho_{\varepsilon}} ( z)$.
\end{enumerate}

\subsection{La $L^2$-estimation locale et r{\'e}solution locale du
$\bar{\partial}$ avec esti- \ \ \ \ \ \ mation jusqu'au bord}

Soit $\omega_0 = i \displaystyle\sum_{j = 1}^n d z_j \wedge d
   \bar{z_{j}} $ la m{\'e}trique k$\mathsf{\ddot{\mathrm{a}}}$hl{\'e}rienne
    standard sur $\mathbb{C}^n $.
Il existe des m{\'e}triques k$\mathsf{\ddot{\mathrm{a}}}$hl{\'e}riennes
$\omega_{\varepsilon}$ sur $U_{\varepsilon}$ , $0 \leqslant \varepsilon
\leqslant \varepsilon_0$, avec les propri{\'e}t{\'e}s suivantes ( lemme 3.1
dans J.Brinkschulte \cite{Jud1}):
\begin{enumerate}
  \item Soit $\gamma_1^{\varepsilon} \leqslant \ldots . \leqslant
  \gamma_n^{\varepsilon}$ les valeurs propres de $i \partial
  \overline{\partial} - \log ( - \rho_{\varepsilon})$ relativement {\`a}
  $\omega_{\varepsilon}$. Alors on a pour $0 < \varepsilon \leqslant
  \varepsilon_0$ , $\gamma_1^{\varepsilon} + \ldots . + \gamma_n^{\varepsilon}
  \geqslant \frac{1}{- \rho_{\varepsilon}}$ .
  
  \item Il existe une constante $M > 0$ telle que \ $\varepsilon
  \frac{\gamma}{q} \omega_0 \leqslant \omega_{\varepsilon} \leqslant M
  \omega_0$ , pour tout \ $0 < \varepsilon \leqslant \varepsilon_0$ . 
\end{enumerate}

Le th{\'e}or{\`e}me suivant nous donne une r{\'e}solution locale du
$\bar{\partial}$ jusqu'au bord dans un domaine $q$-convexe de $\mathbb{C}^n$
(Voir J.Brinkschulte \cite{Jud1}).

\begin{theorem}\label{th5}
  
\end{theorem}

Soit $\Omega \subset \mathbb{C}^n$ un ouvert de classe $\mathcal{C}^2$ , et
supposons que $\Omega$ est $q$-convexe au voisinage de $z_0 \in \partial
\Omega$. Alors il existe un syst{\`e}me fondamental de voisinages $U$ de $z_0$
tel que pour $f \in \mathcal{C}_r^{\infty} ( \overline{\Omega \cap U})$ avec
$\bar{\partial} f = 0$ , $r \geqslant q$ , il existe $u \in \mathcal{C}_{r -
1}^{\infty} ( \overline{\Omega \cap U})$ telle que $\bar{\partial} u = f$.

\begin{lemme}\label{l5}
  (Voir J.Brinkschulte \cite{Jud1})
\end{lemme}

Pour tout $k \in \mathbb{N}$ et $f \in \mathcal{C}_r^{\infty} (
\bar{U}_{\varepsilon k})$ avec $\bar{\partial} f = 0$ , $r \geqslant q$ , il
existe $u \in \mathcal{C}_{r - 1}^{\infty} ( U_{\varepsilon k})$ telle que \
$\bar{\partial} u = f$ et
\[ \forall s \in \mathbb{N} , \text{} \| u \|_{( s + 1,
   U_{\varepsilon_{k + 1}})} \leqslant \frac{M_s}{\varepsilon_k^{s + 2 + r /
   2}} \| f \|_{( s, U_{\varepsilon_k})}, \]
o{\`u} $M_s$ est une constante ind{\'e}pendante de $k$.

\subsection{Preuve du th{\'e}or{\`e}me \ref{th2}}

Dans la preuve du th{\'e}or{\`e}me \ref{th2} on remplacera le r{\'e}sultat de Kohn \cite{Kohn} de r{\'e}solution jusqu'au bord par celui de J.Brinkschulte \cite{Jud1} du
th{\'e}or{\`e}me \ref{th5} qui donne la r{\'e}solution locale jusqu'au bord dans un
domaine $q$-convexe de $\mathbb{C}^n$.\\
\\
\noindent Ecrivons les d{\'e}tails de la preuve :

Soit $f \in W^r ( \Gamma \cap \bar{U})$ \ telle que $\bar{\partial} f = 0$
avec $r \geqslant q$. Puisque la solution est locale, alors on peut prendre
$\Gamma$ born{\'e} d'o{\`u} il existe $R > 0$ tel que $\Gamma = \Gamma_R$.
Ainsi si $z_0 \in \Gamma$ alors $\forall \nu$ , $z_0 \in \Omega_{\nu}^R$.
Consid{\`e}rons $( U_{\varepsilon})_{\varepsilon}$ la famille de voisinages de
$z_0$ donn{\'e} par le lemme \ref{l4} et $\tilde{f}$ une extension de $f$ {\`a}
$\bar{U}_{\varepsilon_1}$ , alors $\bar{\partial} \tilde{f}$ est une $( r +
1)$-forme diff{\'e}rentielle $\bar{\partial}$-ferm{\'e}e dans
$\bar{U}_{\varepsilon_1}$. D{\'e}signons par $h_{\nu}$ une solution de
l'{\'e}quation $\bar{\partial} h_{\nu} = \bar{\partial} \tilde{f}$ dans
$U_{\varepsilon_{\nu}}$ , $\nu \geqslant 1$ , fournie par le lemme \ref{l5} ;

on a donc
\[ \| h_{\nu} \|_{( s + 1, U_{\varepsilon_{\nu + 1}})} \leqslant
   \frac{M_s}{\varepsilon_{\nu}^{s + 2 + r / 2}} \| \bar{\partial} \tilde{f}
   \|_{( s, U_{\varepsilon_{\nu}})} . \]
Consid{\`e}rons sur $U_{\varepsilon_{\nu + 2}}$ la forme diff{\'e}rentielle
$h_{\nu} - \tilde{f}$ dans $U_{\varepsilon_{\nu + 2}}$, qui est
$\bar{\partial}$-ferm{\'e}e ; d{\'e}signons par $\alpha_{\nu}$ une solution de
l'{\'e}quation $\bar{\partial} \alpha_{\nu} = h_{\nu} - \tilde{f}$ fournie par
le lemme \ref{l5}. On a

$\| \alpha_{\nu} \|_{( s + 2, U_{\varepsilon_{\nu + 3}})} \leqslant
\frac{M_{s + 1}}{\varepsilon_{\nu}^{s + 3 + ( r - 1) / 2}} \| h_{\nu} -
\tilde{f} \|_{( s + 1, U_{\varepsilon_{\nu + 2}})}$.

Soit encore
\[ \| \alpha_{\nu} \|_{( s + 2, U_{\varepsilon \nu + 3})} \leqslant \frac{M_s
   M_{s + 1}}{\varepsilon_{\nu}^{2 s + r + 4}}  \| \bar{\partial} \tilde{f}
   \|_{( s, U_{\varepsilon_{\nu}})} . \]
On peut choisir $\tilde{f}$ de sorte que pour tout $s, N \in \mathbb{N}$, on
ait une constante $C$ ind{\'e}pendante de $\nu$ telle que $\| \bar{\partial}
\tilde{f} \|_{( s, U_{\varepsilon_{\nu}})} \leqslant C \| \bar{\partial} f
\|_{( s, U_{\varepsilon_1})} .$ Par cons{\'e}quent
\[ \| \alpha_{\nu} \|_{( s + 2, U_{\varepsilon \nu + 3})} \leqslant \frac{M_s
   M_{s + 1} C b^N}{\varepsilon_{\nu}^{2 s + r + 4 - N}} = K
   \varepsilon_{\nu}^{N - ( 2 s + r + 4)} . \]
On choisit alors $N > 2 s + r + 4$ pour que cette in{\'e}galit{\'e} , jointe
au lemme de Sobolev fournissent la convergence de $\Sigma \alpha_{\nu}$ dans
$W^{r - 1} ( \Gamma_R \cap \bar{U})$. La somme $\alpha$ de la s{\'e}rie
v{\'e}rifie $\bar{\partial} \alpha = f$.

\section{Des exemples de ferm{\'e}s de $\mathbb{C}^n$ v{\'e}rifiant la condi-
\ \ \ \ \ tion $( \lambda')$.}

\begin{itemize}
  \item Les ferm{\'e}s de $\mathbb{C}^n$ v{\'e}rifiant la condition $(
  \lambda)$ de Dufresnoy \cite{Duf} v{\'e}rifient la condition $(
  \lambda')$ car tout domaine pseudoconvexe est $q$-convexe , $1 \leqslant q
  \leqslant n - 1$.
  
  \item L'intersection localement finie d'adh{\'e}rences d'ouverts strictement
  $q$-convexes de $\mathbb{C}^n$ v{\'e}rifie la condition $( \lambda')$. \
  
  En effet on proc{\'e}dera comme dans \cite{Duf}. Ainsi d{\'e}signons par
  $\rho_1, \ldots , \rho_n \ldots$ les fonctions d{\'e}finissant les
  ouverts strictement $q$-convexes $U_n$ de classe $\mathcal{C}^2$ ; autrement
  dit, $\rho_n$ est de classe $\mathcal{C}^2$ dans un voisinage $V_n$ de
  $\mathbb{C}^q$ , le gradiant de $\rho_n$ ne s'annule pas sur $\partial V_n$
  et la restriction au plan tangent complexe de la forme de Levi de $\rho_n$
  admet $q + 1$ valeurs propres strictement positives.
  
  Pour $R$ fix{\'e} , il existe $n_1, \ldots, n_p$ tel que
  \[ \Gamma_R = \{ z \in \mathbb{C}^n  ; \rho_{n_j} ( z) \leqslant 0
     ; | z |^2 - R^2 \leqslant 0 \} . \]
  Quitte {\`a} composer les fonctions $\rho_{n_j}$ par une fonction convexe ,
  on peut supposer que les fonctions $\rho_{n_j}$ sont strictement
  $q$-convexes dans un voisinage de la fronti{\`e}re de $U_{n_j}$ ; donc
  \[ \Omega_{\nu}^R = \left\{ z \in \mathbb{C}^n ; \sup_j \rho_j ( z) <
     \frac{1}{2^{\nu}} \right\} \]
  est un ouvert $q$-convexe pour $\nu$ assez grand et que $\cap \Omega_{\nu}^R
  = \Gamma_R$. \
  
  Du fait que l'enveloppe sup{\'e}rieure d'une famille finie de fonctions
  lipschitziennes est lipschitzienne , il existe une constante $a$ telle que
  \[ \Omega_{\nu}^R \supset \left\{ z \in \mathbb{C}^n ; \rm{dist} ( z,
     \Gamma_R) < \frac{a}{2^{\nu}} \right\} \]
  et la condition sur le gradiant des fonctions $\rho_j$ fournit une constante
  $b$ telle que
  \[ \Omega_{\nu}^R \subset \left\{ z \in \mathbb{C}^n ; \rm{dist} ( z,
     \Gamma_R) < \frac{b}{2^{\nu}} \right\} . \]
  Si $R' > R$ , il existe un voisinage de $\Gamma_{R'}$ tel que
  $\Omega_{\nu}^R$ soit d{\'e}fini dans ce voisinage de $\Gamma_{R'}$ par des
  fonctions $q$-convexes. Si on choisit le voisinage de $\Gamma_{R'}$
  $q$-convexe , on en d{\'e}duit que toute $( q - 1)$-forme
  $\bar{\partial}$-ferm{\'e}e d{\'e}finie dans $\overline{\Omega_{\nu}^R}$ est
  limite uniforme de $( q - 1)$-formes $\bar{\partial}$-ferm{\'e}es dans ce
  voisinage de $\Gamma_{R'}$.
\end{itemize}

\end{document}